# Infinite Divisibility and Max-Infinite Divisibility
# with Random Sample Size

**S. Satheesh**

Department of Statistics
Cochin University of Science and Technology
Cochin 682 022, India.

and

**E. Sandhya**

Prajyoti Niketan College
Pudukkad, Thrissur 680 301, India.

**Abstract:** Infinitely divisible (ID) and max-infinitely divisible (MID) laws are studied when the sample size is random. $\varphi$-ID and $\varphi$-MID laws introduced and studied here approximate random sums and random maximums. The main contributions in this study are: (i) in discussing a class of probability generating functions for $N$, the sample size, (ii) a NaS condition that implies the convergence to an ID (MID) law by the convergence to a $\varphi$-ID ($\varphi$-MID) law and vise versa and thus a discussion of attraction and partial attraction for $\varphi$-ID and $\varphi$-MID laws.

## 1. Introduction

Let $\{X_i\}$ be a sequence of independent and identically distributed (i.i.d) random variables (r.v) and $N$ be a non-negative integer valued r.v independent of $\{X_i\}$. Then $X_1 + \ldots + X_N$ defines a random sum (N-sum) and $\vee\{X_1, \ldots, X_N\}$ defines a random maximum (N-max). Random sums and random maximums have applications in Insurance and Hydrology, (Kaufman (2001)), Finance (Rachev (1993) and Mittnick and Rachev (1993)), Queuing (Makowski (2001)), Reliability (Cai and Kalashnikov (2000)), and Actuaries (Denuit, *et al* (2002)). In the present paper we generalize the facts that the limit distributions of compound Poisson laws are infinitely divisible (ID) and the limit laws of Poisson-maximums are max-infinitely divisible (MID) (Theorem.3.1) to define

ssatheesh@sancharnet.in

φ-ID and φ-MID laws, where $\varphi$ is a Laplace transform (LT), and these laws approximate N-sums and N-maximums.

**Definition 1.1:** A characteristic function (CF) $f(t)$ is φ-ID if for every $\theta \in \Theta$, there exists a CF $h_\theta(t)$, a probability generating function (PGF) $P_\theta$ that is independent of $h_\theta$, such that $P_\theta\{h_\theta(t)\} \rightarrow f(t) = \varphi\{-\log \omega(t)\} \ \forall \ t \in \mathbf{R}$, as θ ↓ 0 through a $\{\theta_n\} \in \Theta$. Here, $\varphi$ is the LT of a r.v $Z > 0$ and $\omega(t)$ is a CF that is ID.

**Definition 1.2:** A distribution functions (d.f) $F$ on $\mathbf{R}^d$ for $d \geq 2$ integer, is φ-MID if for every $\theta \in \Theta$, there exists a d.f $H_\theta$, a PGF $P_\theta$ that is independent of $H_\theta$, such that $P_\theta\{H_\theta\} \rightarrow F = \varphi\{-\log G\}$ as $\theta \downarrow 0$ through a $\{\theta_n\}$. Here, $\varphi$ is the LT of a r.v $Z > 0$ and $G$ is a d.f that is MID.

This study is a continuation of Satheesh (2001*a*,*b* & 2002). Generalizations of ID and geometrically-ID laws were considered by: Sandhya (1991, 1996), discussing two examples of non-geometric laws for $N$ but the description was not constructive. The descriptions of N-ID laws (with CF $\varphi\{\psi\}$, where exp$\{-\psi\}$ is a CF that is ID) by Gnedenko and Korolev (1996), Klebanov and Rachev (1996) and Bunge (1996), are based on the assumption that the PGF $\{P_\theta, \theta \in \Theta\}$ of $N_\theta$ formed a commutative semi-group. Satheesh (2001*b*) showed that this assumption is not natural since it captures the notion of stable rather than ID laws. Also, it rules out any $P_\theta$ having an atom at the origin. A glaring situation is that this theory cannot approximate negative binomial sums by taking $\varphi$ as gamma as one expects. Kozubowski and Panorska (1996, 1998) discussed ν-stable laws with CF $\varphi\{\psi\}$ (exp$\{-\psi\}$ being a stable CF) approximating $N_\theta$-sums, assuming $N_\theta \xrightarrow{\ p\ } \infty$ as $\theta \downarrow 0$ and handled negative binomial sums. But identifying those $N_\theta \xrightarrow{\ p\ } \infty$ was not done. Satheesh, *et al* (2002) showed that when $\varphi$ is gamma, $N_\theta$ is Harris with PGF $s/[a - (a-1) \ s^k]^{1/k}$, $k > 0$ integer, $a > 1$, for its N-sum stability. Satheesh (2001*b*), motivated by these observations and to overcome them introduced φ-ID laws approximating $N_\theta$-sums where the PGF of $N_\theta$ is a member of $\wp_\varphi(s) = \{P_\theta(s) = s^j \ \varphi \ [(1 - s^k)/\theta] \}$, $0 < s < 1$, $j \geq 0$ & $k \geq 1$ integer and $\theta > 0\}$ without requiring that $N_\theta \xrightarrow{\ p\ } \infty$ and discussed φ-attraction and partial φ-attraction. He also proved an analogue of Theorem 4.6.5 of Gnedenko and Korolev (1996, p.149) where the convergence to an ID law implies the convergence to a φ-ID law and vise versa for the class of PGFs $\wp_\varphi(s)$. The gain is: here for every LT $\varphi$ we have a class of PGFs w.r.t which we can discuss the φ-ID law unlike those in the previous works where either the PGF is unique for a LT (the commutative semi-group approach) or we need the assumption $N_\theta \xrightarrow{\ \ p\ \ } \infty$ as $\theta \downarrow 0$.

In the literature we also have the discussion of MID laws parallel to the ID laws. MID laws were studied by Balkema and Resnick (1977), geometrically MID laws by Rachev and Resnick (1991), Mohan (1998), and Satheesh (2002) discussed φ-ID and φ-MID laws from the perspective of randomization and mixtures of ID and MID laws. Also see the G-max-stable laws by Sreehari (1995), and the stability of N-extremes of continuous and discrete laws by Satheesh and Nair (2002).

To improve upon these results, in Section.2 we discuss φ-ID laws giving a necessary and sufficient condition for a CF to be φ-ID w.r.t the class of PGFs $\wp_\varphi(s)$ of Satheesh (2001$b$) for $N$, the sample size. In Section.3 we describe φ-MID laws paralleling φ-ID laws that generalizes the geometrically MID laws of Rachev and Resnick (1991), Mohan (1998) and the G-max-stable laws of Sreehari (1995).

## 2. Ramifications of Feller's Proof of Bernstein's Theorem and φ-ID laws

The following two lemmas are ramifications of Feller's proof of Bernstein's Theorem (Feller (1971) p.440) and were first proved by Satheesh (2001$b$). Together these two lemmas provide a class of discrete laws that can be used in the transfer theorems for random sums and maximums of Gnedenko (1982). Notice that these theorems hold true for random vectors on $\boldsymbol{R}^d$ for $d \geq 2$ integer, as well.

**Lemma 2.1:** $\wp_\varphi(s) = \{P_\theta(s) = s^j \varphi\{(1-s^k)/\theta \}, 0<s<1, j \geq 0$ & $k \geq 1$ integer and $\theta >0\}$ describes a class of PGFs for any given LT $\varphi$.

**Proof:** Since $P_\theta(s)$ is absolutely monotone and $P_\theta(1) = 1$, by Feller (1971, p.223) $P_\theta(s)$ is a PGF.

**Lemma 2.2:** Given a r.v $U$ with LT $\varphi$, the integer valued r.vs $N_\theta$ with PGF $P_\theta$ in the class $\wp_\varphi(s)$ described in Lemma.2.1 satisfy

$\theta N_\theta \xrightarrow{d} kU$ as $\theta \to 0$.

**Proof:** The LT of $\theta N_\theta$ is $e^{-vj\theta}\varphi\{(1-e^{-vk\theta})/\theta\}$, $v > 0$ and by Feller (1971, p.440),

$$\underset{\theta \to 0}{Lt}\{e^{-vj\theta}\varphi((1-e^{-vk\theta})/\theta)\} = \varphi(kv).$$

We now observe that perhaps the assumption, $\{P_\theta : \theta \in \Theta\}$ form a commutative semi-group, is not a natural setting for the description of N-ID laws.

**Remark.2.1:** Suppose the $N_\theta$-sum of a LT $\varphi(s)$ is of the same type. That is $\varphi(s) = P_\theta\{\varphi(\theta s)\}$. Equivalently, $\varphi(s/\theta) = P_\theta\{\varphi(s)\}$, and when $s = \varphi^{-1}(z)$ we have

$$P_\theta(z) = \varphi\{\varphi^{-1}(z)/\theta\}, \; \theta \in \Theta.$$

This relation thus captures the structure of N-sum-stability. By corollary.4.6.1 in Gnedenko and Korolev (1996, p.141) this relation is equivalent to the assumption that $\{P_\theta(z) : \theta \in \Theta\}$ form a commutative semi-group. Clearly the converse also holds. Notice that Gnedenko & Korolev (1996) and Klebanov & Rachev (1996) first developed N-normal laws and then generalized it to describe N-ID laws. But normal laws are not just ID but stable as well and this is one reason why the commutative semi-group assumption came in to the picture here.

Now recall our description of $\varphi$-ID laws in definition.1.1. Satheesh (2001*b*) showed that a $\varphi$-ID CF $f(t) = \varphi\{-\log \omega(t)\}$ has no real zeroes. Now, let $\{X_{\theta,j} ; \theta \in \Theta, j \geq 1\}$ be a sequence of i.i.d r.vs with CF $h_\theta$ and for $k \geq 1$ integer, set $S_{\theta,k} = X_{\theta,1} + \ldots + X_{\theta,k}$. In the sequel we will consider sequences $\{\theta_n \in \Theta\}$ and $\{k_n ; n \geq 1\}$ of natural numbers such that as $\theta \downarrow 0$ through $\{\theta_n\}$ and $S_{\theta,k_n} \xrightarrow{d} U$ so that $U$ is ID. Also in an $N_\theta$-sum of $\{X_{\theta,k}\}$, $N_\theta$ and $X_{\theta,k}$ are assumed to be independent for each $\theta \in \Theta$.

**Example 2.1:** Let $\varphi$ be a LT. Then $\varphi\{(1-s)/\theta\}$, $\theta > 0$ is a PGF and the CF $f(t) = \underset{n \to \infty}{Lt}\varphi\{a_n(1-h_n(t))\}$ is $\varphi$-ID, where $\{a_n\}$ are some positive constants and $\{h_n(t)\}$ are CFs and the distributions of $\{h_n(t)\}$ and $\varphi\{a_n(1-s)\}$ are independent for each $n \geq 1$ integer. Recall the classical de-Finetti theorem here.

**Theorem 2.1:** The limit law of $N_\theta$-sums of $\{X_{\theta,j}\}$ as $\theta \downarrow 0$ through a $\{\theta_n\}$, where the PGF of $N_\theta$ is a member of $\wp_\varphi(s)$, is necessarily φ-ID. Conversely, for any given LT $\varphi$ ; the φ-ID law can be obtained as the limit law of $N_\theta$-sums of i.i.d r.vs as $\theta \downarrow 0$ for each member of $\wp_\varphi(s)$.

**Proof:** Follows from Lemma 2.1, 2.2 and the transfer theorem for sums.

Next we prove a necessary and sufficient condition for $f(t) = \underset{\theta \to 0}{Lt}\, P_\theta\{h_\theta(t)\}$ to be φ-ID in terms of the CF $h_\theta$, when the PGF of $N_\theta \in \wp_\varphi(s)$. This incidentally implies the convergence of non-random sum of $\{X_{\theta,i}\}$ to an ID law and vice-versa as well. From another angle it generalizes Theorem.1.1 of chapter.17 in Feller ((1971) p.555).

**Theorem 2.2:** A CF $f = \underset{\theta \to 0}{Lt}\, \{h_\theta\}^j\, \varphi((1- \{h_\theta\}^k)/\theta )$, $j \geq 0$ & $k \geq 1$ integer, is φ-ID iff there exists a continuous function $\psi(t)$ such that $\forall\, t \in \boldsymbol{R}$

$$(1- h_\theta(t))/\theta \;\to \psi(t), \text{ as } \theta \downarrow 0 \text{ through a sequence } \{\theta_n\}.$$

In this case $f = \varphi(\psi)$ .

**Proof:** The condition is sufficient since it implies $\underset{\theta \to 0}{Lt}\, h_\theta = 1$ and that $e^{-\psi}$ is ID by Feller (1971, p.555) and invoking Lemma.2.2 and the transfer theorem for sums.

To prove the necessity we have: $\underset{\theta \to 0}{Lt}\, \{h_\theta\}^j\, \varphi((1- \{h_\theta\}^k)/\theta ) = f(t)$ is φ-ID. Then there exists a CF $\omega$ that is ID (hence no real zeroes), such that $-\log\{\omega\} = \psi$ and consider the CF $\omega^k$, $k \geq 1$ integer that is again ID. Setting $f = \varphi(k\psi) = \varphi(-k \log\{\omega\})$;

$$\underset{\theta \to 0}{Lt}\, \{h_\theta\}^j\, \underset{\theta \to 0}{Lt}\, \varphi(\tfrac{1}{\theta}(1- \{h_\theta\}^k)) = \underset{\theta \to 0}{Lt}\, \{h_\theta\}^j\, \varphi(\underset{\theta \to 0}{Lt}\, \tfrac{1}{\theta}(1- \{h_\theta\}^k)) = \varphi(-k \log\{\omega\}).$$

Now, since $\theta \downarrow 0$ the existence of $\log\{\omega\}$ on the RHS is guaranteed only if $\underset{\theta \to 0}{Lt}\, (1- \{h_\theta\}^k)$ $= 0$. Hence $\underset{\theta \to 0}{Lt}\, h_\theta^k = 1$ so that $\underset{\theta \to 0}{Lt}\, h_\theta = 1$. Hence remembering that $\log(1-z) \sim -z$ as $z \downarrow 0$ we get:

$$\underset{\theta \to 0}{Lt} \; \tfrac{1}{\theta}(1-\{h_\theta\}^k) \; = -k\log\{\omega\} \qquad \Rightarrow \qquad \underset{\theta \to 0}{Lt} \; \tfrac{1}{\theta}\log\{h_\theta\}^k = k\log\{\omega\}$$

$$\Rightarrow \qquad \underset{\theta \to 0}{Lt} \; \tfrac{k}{\theta}\log\{h_\theta\} = k\log\{\omega\} \qquad \Rightarrow \qquad \underset{\theta \to 0}{Lt} \; \tfrac{1}{\theta}\log\{h_\theta\} = \log\{\omega\} \qquad (*)$$

which again imply:  $\underset{\theta \to 0}{Lt} \log\{h_\theta\} = 0 \qquad \Rightarrow \qquad \underset{\theta \to 0}{Lt} \; h_\theta = 1.$

Hence : $\tfrac{1}{\theta}\log\{h_\theta\} \; = \; -\tfrac{1}{\theta}(1-\{h_\theta\}) \, [1+o(1)],$

where $o(1)$ is a quantity that vanishes as $\theta \downarrow 0$. By $(*)$ the LHS tends to $\log\{\omega\}$ and hence

$$\underset{\theta \to 0}{Lt} \; (1-\{h_\theta\})/\theta \; = -\log\{\omega\} \; = \psi \; , \text{ as was to be proved.}$$

In the classical summation scheme a CF $h(t)$ belongs to the domain of attraction (DA) of the CF $\omega(t)$ if there exists a sequence of real constants $a_n > 0$ and $b_n$ such that

$$\text{as } \; n \to \infty \, , \exp\{-n(1-h_n(t))\} \to \exp\{-\psi(t)\} = \omega(t) \; \forall \; t \in \boldsymbol{R}.$$

where $h_n(t) = h(t/a_n)\exp(-itb_n)$. A CF $h(t)$ belongs to the domain of partial attraction (DPA) of the CF $\omega(t)$ if there exists an increasing subsequence $\{n_m\}$ of positive integers such that

$$\text{as } m \to \infty, \exp\{-n_m(1-h_m(t))\} \to \exp\{-\psi(t)\} = \omega(t) \; \forall \; t \in \boldsymbol{R}.$$

Now, we generalize these notions w.r.t PGFs that are assumed to be independent of the CF $h(t)$, the r.vs $N_\theta$ corresponding to these PGFs $P_\theta$ satisfying $\theta N_\theta \xrightarrow{d} U$ as $\theta \downarrow 0$ and $\varphi$ being the LT of the r.v $U$. We formulate this by considering $\{\theta_n \in \Theta\}$ such that as $n \to \infty$, $\theta \downarrow 0$ through $\{\theta_n\}$. The PGFs will be denoted by $P_n = P_{\theta_n}$ corresponding to $\theta_n = 1/n$, where $\{n\}$ is the sequence of positive integers and $P_{n_m} = P_{\theta_m}$ corresponding to $\theta_m = 1/n_m$, a subsequence $\{n_m\}$ of $\{n\}$. Now with $h_n(t)$ as defined above;

**Definition 2.1:** A CF $h(t)$ belongs to the domain of φ-attraction (Dφ-A) of the CF $f(t)$ if there exists sequences of real constants $a_n = a(\theta_n) > 0$ & $b_n = b(\theta_n)$ such that $\underset{n \to \infty}{Lt} P_n\{h_n(t)\} = f(t) \ \forall \ t \in \mathbf{R}$ and $h(t)$ is in the domain of partial φ-attraction (DPφ-A) of $f(t)$ if; $\underset{m \to \infty}{Lt} P_{n_m}\{h_m(t)\} = f(t) \ \forall \ t \in \mathbf{R}$.

Certain implications of Theorem 2.2 are the following: This theorem enables us to conclude (w.r.t the PGFs in $\wp_\varphi(s)$) that if the CF $h(t)$ belongs to the DA of the stable law (DPA of the ID law) with CF $\omega(t) = e^{-\psi(t)}$, then it is also a member of Dφ-A of the φ-stable law (DPφ-A of a φ-ID law) with CF $f(t) = \varphi\{\psi(t)\}$ and the converses are also true. All that we need is to prescribe, $\theta_n = 1/n$ for φ-attraction and $\theta_m = 1/n_m$ for partial φ-attraction. Thus the DA of a stable law (DPA of an ID law) with CF $\omega(t) = e^{-\psi(t)}$ coincides with the Dφ-A of the φ-stable law (DPφ-A of a φ-ID law) with CF $f(t) = \varphi\{\psi(t)\}$ for each PGF $P_\theta \in \wp_\varphi(s)$, and none of them are empty as well. Satheesh (2001$b$) also has conceived these notions and reached these conclusions, first for the class of PGFs $\varphi\{(1-s)/\theta\}$, $\theta > 0$ and then w.r.t the PGFs in $\wp_\varphi(s)$ by proving an analogue of Theorem 4.6.5 of Gnedenko and Korolev ((1996) p.149). These notions and conclusions generalize those on attraction and partial attraction in geometric sums in Sandhya (1991), Sandhya and Pillai (1999), Mohan *et al*. (1993), Ramachandran (1997), and those for N-sums in Gnedenko and Korolev (1996). One may extend these results to random vectors on $\mathbf{R}^d$ for $d \geq 2$ integer, generalizing proposition 2.1 of Kozubowski and Panorska (1998). Satheesh (2003) has discussed operator φ-stable laws and its advantages over operator ν-stable laws in its application to continuous time random walks.

### 3. φ-MID laws

Since all d.fs on $\mathbf{R}$ are max infinitely divisible (MID) discussion of MID laws and their generalizations are relevant only for d.fs on $\mathbf{R}^d$ for $d \geq 2$ integer. Further, the operations are to be taken component wise. By $\{F > 0\}$ we denote the set $\{y \in \mathbf{R}^d : F(y) > 0\}$ and the discussion is w.r.t $\{F > 0\}$. From Balkema and Resnick (1977) we have: $G$ is MID iff $G$ has an exponent measure $\mu$ and $\{G > 0\}$ is a rectangle. Also let $\lambda = $ inf $\{G > 0\}$ be the bottom of this rectangle. Notice that the transfer theorem for maximums holds true for random vectors also. From Balkema and Resnick (1977, Theorem.3 and Corollaries 3 & 4) we have:

**Theorem 3.1:** $G$ is MID iff $G = \underset{n \to \infty}{Lt} \exp\{-a_n(1 - H_n)\}$, for some d.fs $\{H_n\}$ and constants $\{a_n > 0\}$.

**Theorem 3.2:** $G$ is MID iff there exists a $\lambda \in [-\infty, \infty)^d$ and an exponent measure $\mu$ concentrating in $[\lambda, \infty]$ such that for $y \geq \lambda$ $G(y) = \exp\{-\mu([\lambda, \infty]^c)\}$.

Now recall our description of $\varphi$-MID laws in definition.1.2.

**Theorem 3.3:** The limit of a sequence of $\varphi$-MID laws is again $\varphi$-MID.

**Proof:** Follows from the following relation where $\{F_n\}$ are $\varphi$-MID, $\{G_n\}$ are MID so that $G$ is MID;

$$F = \underset{n \to \infty}{Lt} F_n = \underset{n \to \infty}{Lt} \varphi\{-\log G_n\} = \varphi\{-\log \underset{n \to \infty}{Lt} G_n\} = \varphi\{-\log G\}.$$

**Theorem 3.4:** For a d.f $F$ on $\boldsymbol{R}^d$, the following are equivalent.

(i) $F$ is $\varphi$-MID

(ii) $\exp\{-\varphi^{-1}(F)\}$ is MID

(iii) There exists a $\lambda \in [-\infty, \infty)^d$ and an exponent measure $\mu$ concentrating in $[\lambda, \infty]$ such that for $y \geq \lambda$, $F(y) = \varphi\{\mu([\lambda, y]^c\}$

(iv) There exists a multivariate extremal process $\{Y(t), t > 0\}$ governed by a MID law and an independent r.v $Z$ with d.f $F_z$ and LT $\varphi$ such that $F(y) = P\{Y(Z) \leq y\}$.

**Proof:** (i) $\Rightarrow$ (ii) $F$ is $\varphi$-MID implies $F = \varphi\{-\log G\}$ where $G$ is MID $\Rightarrow \exp\{-\varphi^{-1}(F)\}$ = $G$ is MID.

(ii) $\Rightarrow$ (iii) By (ii) $\exp\{-\varphi^{-1}(F)\} = G = \exp\{-\mu([\lambda, y]^c)\}$. Hence; $\varphi^{-1}(F) = \mu([\lambda, y]^c)$ or $F(y) = \varphi\{\mu([\lambda, y]^c\}$.

(iii) $\Rightarrow$ (iv) By (iii) we have an exponent measure $\mu$ and let $\{Y(t), t>0\}$ be the extremal process governed by a MID law with exponent measure $\mu$. That is; $P\{Y(t) \leq y\}$ = $\exp\{-t\mu([\lambda,y]^c)\}$. Hence

$$P\{Y(Z) \leq y\} = \int_{o}^{\infty} \exp\{-t\mu([\lambda,y]^c)\} \, dF_Z(t) = \varphi\{\mu([\lambda,y]^c\} = F(y).$$

(iv) $\Rightarrow$ (i) This is now obvious. Thus the proof is complete.

This result generalizes Proposition.2.2 in Rachev and Resnick (1991) when $N_\theta$ is geometric. Analogous to their Proposition.3.2 we can have a result for $\varphi$-stable laws as well.

Now, let $\{Y_{\theta,j} ; \theta \in \Theta, j \geq 1\}$ be a sequence of i.i.d random vectors in $\boldsymbol{R}^d$, $d \geq 2$ integer with d.f $H_\theta$ and for $k \geq 1$ integer, set $\boldsymbol{M}_{\theta,k} = \vee \{Y_{\theta,1} , \ldots, Y_{\theta,k}\}$. Again we will consider a $\{\theta_n \in \Theta\}$ and $\{k_n ; n \geq 1\}$ of natural numbers such that as $\theta \downarrow 0$ through $\{\theta_n\}$ and $\boldsymbol{M}_{\theta, k_n} \xrightarrow{d} V$ so that $V$ is MID. Also in an $N_\theta$-max of $\{Y_{\theta,j}\}$, $N_\theta$ and $Y_{\theta,j}$ are independent for each $\theta \in \Theta$.

**Remark 3.1:** When the integer valued r.vs $N_\theta$ with PGF $P_\theta$ satisfies $\theta N_\theta \xrightarrow{d} Z$ as $\theta \downarrow 0$ and $\underset{n\to\infty}{Lt} \{H_\theta\}^n = G$, $P_\theta\{H_\theta\} \to F = \varphi\{-\log G\}$ is always satisfied by invoking the transfer theorem for maximums. Because by Balkema and Resnick (1977) if $\{H_n\}^n \to G$, then $G$ is MID.

**Example 3.1:** Let $\varphi$ be a LT. Then $\varphi\{(1-s)/\theta\}$, $\theta >0$ is a PGF and the d.f $F = \underset{n\to\infty}{Lt} \varphi\{a_n(1- H_n)\}$ is $\varphi$-MID, where $\{a_n\}$ are some positive constants and $\{H_n\}$ are d.fs and the distributions of $\{H_n\}$ and $\varphi\{a_n(1- s)\}$ are independent for each $n$. Recall Theorem.3.1 here.

**Theorem 3.5:** $F$ is the limit law of $N_\theta$-maxs of i.i.d random vectors as $\theta \downarrow 0$, where the PGF of $N_\theta$ is $\varphi\{(1-s)/\theta\}$, iff $F$ is $\varphi$-MID.

**Theorem 3.6:** The limit law of $N_\theta$-maxs of i.i.d random vectors as $\theta \downarrow 0$, where the PGF of $N_\theta$ is a member of $\wp_\varphi(s)$, is necessarily $\varphi$-MID. Conversely, for any given LT $\varphi$; the $\varphi$-MID law is the limit law of $N_\theta$-maxs of i.i.d random vectors as $\theta \downarrow 0$ for each member of $\wp_\varphi(s)$.

These theorems follow from Lemma 2.1, 2.2 and the transfer theorem for maximums.

**Theorem 3.7:** A d.f $F = \underset{\theta \to 0}{Lt} \{H_\theta\}^j \, \varphi((1- \{H_\theta\}^k)/\theta)$, $j \geq 0$ & $k \geq 1$ integer, is $\varphi$-MID iff there exists a $\lambda \in [-\infty, \infty)^d$ and an exponent measure $\mu$ concentrating in $[\lambda, \infty]$ such that for $y \geq \lambda$

$$(1- H_\theta)/\theta \to \mu, \text{ as } \theta \downarrow 0 \text{ through a sequence } \{\theta_n\}.$$

In this case $F = \varphi(\mu)$.

**Proof:** The condition is sufficient since it implies $\underset{\theta \to 0}{Lt} H_\theta = 1$ and that $e^{-\mu}$ is MID by theorems 3.1 and 3.2 and invoking lemma.2.2 and the transfer theorem for maximums.

To prove the necessity we have: $\underset{\theta \to 0}{Lt} \{H_\theta\}^j \, \varphi((1- \{H_\theta\}^k)/\theta) = F$ is $\varphi$-MID. Then there exists a d.f $G$ that is MID with exponent measure $\mu$, such that $-\log\{G\} = \mu$ and consider the d.f $G^k$, $k \geq 1$ integer that is again MID. Setting

$$F = \varphi(k\mu) = \varphi(-k \log\{G\});$$

$$\underset{\theta \to 0}{Lt} \{H_\theta\}^j \underset{\theta \to 0}{Lt} \varphi(\tfrac{1}{\theta}(1- \{H_\theta\}^k)) = \underset{\theta \to 0}{Lt} \{H_\theta\}^j \, \varphi(\underset{\theta \to 0}{Lt} \tfrac{1}{\theta}(1- \{H_\theta\}^k)) = \varphi(-k \log\{G\}).$$

Now since $\theta \downarrow 0$ the existence of $\log\{G\}$ on the RHS is guaranteed only if $\underset{\theta \to 0}{Lt}(1- \{H_\theta\}^k)$ $= 0$. Hence $\underset{\theta \to 0}{Lt} H_\theta^k = 1$ so that $\underset{\theta \to 0}{Lt} H_\theta = 1$. Hence, as in the proof of Theorem 2.2:

$$\underset{\theta \to 0}{Lt} \tfrac{1}{\theta}(1- \{H_\theta\}^k) = -k \log\{G\} \quad \Rightarrow \quad \underset{\theta \to 0}{Lt} \tfrac{1}{\theta}\log\{H_\theta\} = \log\{G\}. \qquad (*)$$

Hence we have: $\frac{1}{\theta}\log\{H_\theta\} = -\frac{1}{\theta}(1-\{H_\theta\})\,[1+o(1)]$,

where $o(1)$ is a quantity that vanishes as $\theta \downarrow 0$. By $(*)$ the LHS tends to $\log\{G\}$ and hence

$$\underset{\theta\to 0}{Lt}\,(1-\{H_\theta\})/\theta = -\log\{G\} = \mu \;,\; \text{as was to be proved.}$$

The notions of attraction and geometric attraction for maximums have been discussed in Rachev and Resnick (1991) and Mohan (1998). Now we can develop the notions of domain of φ-max-attraction (Dφ-MA) and domain of partial φ-max-attraction (DPφ-MA) for $N_\theta$-maximums paralleling those in the case of $N_\theta$-sums. Again we consider $\{\theta_n \in \Theta\}$ such that as $n\to\infty$, $\theta \downarrow 0$ through $\{\theta_n\}$ with the PGF $P_n$ corresponding to $\theta_n = 1/n$ and $P_{n_m}$ to $\theta_m = 1/n_m$ as before.

**Definition 3.1:** A d.f $H$ in $\boldsymbol{R}^d$, $d \geq 2$ integer, belongs to the Dφ-MA of the d.f $F$ with non-degenerate marginal distributions if there exists normalizing constants $a_{i,n} = a_i(\theta_n) > 0$ and $b_{i,n} = b_i(\theta_n)$ such that with $H_n(y) = H(a_{i,n}y_i + b_{i,n}$, $1 \leq i \leq d)$; $\underset{n\to\infty}{Lt}\,P_n\{H_n\} = F$ and $H$ belongs to the DPφ-MA of the d.f $F$ if; $\underset{m\to\infty}{Lt}\,P_{n_m}\{H_m\} = F$.

Thus w.r.t the PGFs in $\wp_\varphi(s)$, (by invoking theorem.3.7) if the d.f $H$ belongs to the DMA of a max-stable law (DPMA of a MID law) with d.f $G = \mathrm{e}^{-\mu}$ then it is also a member of the Dφ-MA of a φ-max-stable law (DPφ-MA of a φ-MID law) with d.f $F = \varphi\{\mu\}$ and the converses are also true. All that we need is to prescribe, $\theta_n = 1/n$ for φ-max-attraction and $\theta_m = 1/n_m$ for partial φ-max-attraction. Thus the DMA of a max-stable law (DPMA of a MID law) with d.f $G = \mathrm{e}^{-\mu}$ coincides with the Dφ-MA of a φ-max-stable law (DPφ-MA of a φ-MID law) with d.f $F = \varphi\{\mu\}$ for each PGF $P_\theta \in \wp_\varphi(s)$, and none of them are empty also. These results generalize those in Rachev and Resnick (1991) and Mohan (1998).

**Acknowledgement**

Authors wish to thank the referee for certain useful suggestions and Prof. N Unnikrishnan Nair, for the encouragement.